%\magnification=\magstep1
%version August 24, 2011

\def\div{\hbox{\rm div}\,}
\def\curl{\hbox{\rm curl}\,}

\def\vf{\varphi}
\def\ve{\varepsilon}

\def\moh{$(-1)\,$-- homogeneous}
\def\rmo#1{R^{#1}\setminus\{0\}}

 \font\eightrm=cmr8

 \def\ul{u_{\lambda}}
 \def\fl{f_{\lambda}}

 \def\s2{S^2}
 \def\sm{S^{n-1}}
 \def\lrn{\Delta}
 \def\lsm{\Delta_{\hbox{\fiverm S}}}
 \def\ric{\rm{Ric}}
 \def\difsm{\nabla^{\hbox{\fiverm S}}}
 \def\difc{\nabla^{\hbox{\fiverm C}}}
 \def\hl{\Delta_{\hbox{\fiverm H}}}
 \def\cl{\Delta_{\hbox{\fiverm C}}}
 \font\eightrm=cmr8
 \font\fiverm=cmr5
 \def\mez{$\,$}
 \def\pul{{1\over 2}}
 \def\intzpp{\int_0^{\pi/2}}
 \def\NSE{1}
 \def\Stn{2}
 \def\S{3}
 \def\Eone{4}
 \def\Etwo{5}
 \def\Ethree{6}
 \def\Efour{7}
 \def\Efive{8}
 \def\SC{9}
 \def\Esix{10}
 \def\Eseven{11}
 \def\Stwo{12}
 \def\C{13}
 \def\NSEdiv{14}
 \def\SCtwod{15}
 \def\Gtwo{16}
 \def\St{17}
 \def\ODEone{18}
 \def\E{19}
 \def\rc{20}
 \def\P{21}
 \def\Cprime{22}
 \def\der{23}
 \def\efi{24}
 \def\B{25}
 \def\SnB{26}
 \def\Aone{27}
 \def\Atwo{28}
 \def\Sn{29}
 \def\EF{30}
 \def\dec{31}

\parindent=0pt
\advance \baselineskip by .7pt
%\font\eightrm=cmr8
 \centerline{\bf ON LANDAU'S
SOLUTIONS OF THE NAVIER-STOKES EQUATIONS}
\smallskip
\bigskip
\centerline{Vladim\'\i{}r \v Sver\'ak}

 \centerline{University of
Minnesota}

\bigskip
\centerline{\sl Dedicated to Professor Nikolai V. Krylov on the occasion of his 70th birthday.}

\bigskip

 {\bf 1. Introduction.}

 \medskip
 In this paper we study a special class
 of solutions of the $n-$dimensional steady-state  Navier-Stokes equations
 $$ \eqalign{-\Delta u +u\nabla u +\nabla p & =0 \,,\cr
            \div u & = 0\,,\cr} \eqno\hbox{\rm(\NSE)}
            $$
 where $u=(u_1,\dots,u_n)$.
The equations have a non-trivial scaling symmetry $u(x)\to\lambda
u(\lambda x)$ and it is natural to try to find solutions which are
invariant under this scaling.
The simplest natural domain of definition for such solutions
is $R^n\setminus\{0\}$. In this case, assuming that the solutions are smooth in $R^n\setminus\{0\}$, we are able to obtain a good classification of the invariant solutions in all dimensions. There are some interesting conclusions for the regularity theory as well as for the long-distance behavior of solutions in exterior domains which can be drawn from this classification, which will be discussed.  We will distinguish three cases, namely $n=2$, $n=3$ and $n\ge 4$. (Sometimes it is also useful to distinguish the cases $n=4$ and $n\ge 5$, as $n=5$ is the lowest dimension  in which the \moh{} functions which are smooth in $R^n\setminus\{0\}$ have locally finite energy $\int_{|x|<r}|\nabla u|^2\,.$ This will not be important for our purposes in this paper, however.)

\medskip
It is useful to note that the effects of a \moh{} singularity are more serious in low dimensions.
Our \moh{} solution $u$ will be locally integrable across the origin for any $n\ge 2$ and hence can always be considered as a distribution in $R^n$, but the validity of the equations across the origin will depend on $n$. For example, for $n\ge3 $ the equation $\div u=0$ will be  satisfied in $R^n$, but this may not be the case for $n=2$, when case $\div u$ may produce a multiple of a Dirac mass at $x=0$. Similarly, for \moh{} smooth Navier-Stokes solutions in $R^n\setminus\{0\}$,
the weak form of the equation ($\int_{R^n} (-u\Delta\vf-u_iu_j\vf_{i,j}) =0$ for smooth compactly supported vector fields $\vf$ with $\div \vf=0$) will be satisfied across the origin when $n\ge 4$.
For $n=3$ the expression may produce a non-trivial right-hand side (a multiple of a Dirac mass supported at $\{0\}$), while for $n=2$ the right-hand side may not  be well-defined, in general.

\medskip
Let us start with the case $n=3$, which is perhaps the most interesting.
Explicit examples of $(-1)$ homogeneous solutions in $R^3\setminus\{0\}$
were first calculated by L.D.Landau in 1944 ([L]) and can be found
in standard textbooks ([LL], p.\ 82, or [B],
p.\ 206, for example). See also formulae (\Eseven) in Section 4. The main idea of Landau's calculation is that if we
impose an additional symmetry requirement, namely that the
solutions are axi-symmetric, the system (\NSE) reduces to a system
of ODEs which, surprisingly, can be solved explicitly in terms of
elementary functions. (In fact, as it was kindly pointed out to the
author by V. Galaktionov, the ODEs were written down already in 1934
by N.A.Slezkin, see [Sl].) The solutions were also independently found
by H.B.Squire in 1951 ([Sq]). More recently, the topic has been
re-visited in [TX] and [CK], where issues concerning Landau's
solutions are addressed from a slightly different viewpoint.

\medskip
Here we prove that even if we drop the requirement of
axi-symmetry, Landau's solutions are still the only solutions of
(\NSE) which are invariant under the natural scaling. More
precisely, we will prove the following:

\medskip
{\bf Theorem 1.} {\sl Assume that $u\colon R^3\setminus \{0\}\to
R^3$ is a non-trivial smooth solution of (\NSE) satisfying $\lambda
u(\lambda x)=u(x)$ for each $\lambda >0$. Then $u$ is a Landau
solution. In other words, $u$ is axi-symmetric and, in a suitable
coordinate frame, is described by  formulae (\Eseven) in Section 4.}
\medskip

The proof of the theorem shows a connection between the
scale-invariant solutions of (\NSE) and the conformal geometry of
the two-dimensional sphere. In fact, once the connection is
understood, the formulae for Landau's solutions can be derived
without much calculation, using just the geometrical properties of
the two-dimensional sphere.

\medskip

Some implications of Theorem 1 are considered in Sections 2 and 3.

\medskip
In the case $n=2$ our assumptions reduce the problem to an ODE on the circle $S^1$. The ODE has been studied in a 1917 paper by G. Hamel [Ha] where a reasonably complete description of solutions is obtained in terms of elliptic functions. Here we re-visit some of these calculations and classify the solutions satisfying an additional constraint that $\div u =0$ across the origin, which means that the origin is neither a source nor a sink for the flow. It turns out that with this assumption there is, modulo rotations, a countable family of \moh{} solutions (smooth away from the origin). See Theorem 2 in Section 5.

\medskip
The case $n\ge 4$ has been previously considered by several authors in connection with potential singularities of Navier-Stokes solutions. In particular, a non-trivial \moh{} solution which is smooth in $\rmo 5$ would represent singular weak solution in $\{|x|<1\}$ with finite energy $\int_{|x|< 1}|\nabla u|^2$. It was proved independently by several authors ([FR], [St2], [T], [Sv2] ) that for $n\ge 4$ there are no non-trivial \moh{} solutions (smooth away from the origin), thus ruling our this particular scenario for singularities. For completeness we reproduce here the proof  given by the author and T.-P. Tsai, which appeared in [T]. See Theorem 3 in Section 6.

\medskip
The problem of finding \moh{} solutions can be considered in any domain  invariant under the dilations $x\to\lambda x$, and the next domain (after $R^n\setminus\{0\}$) which one should consider is the half-space $R^n_+$, with $u$ satisfying the boundary condition $u=0$ on $\partial R^n_+\setminus\{0\}$. While the problem should be manageable for $n=2$, it seems to quite harder when $n\ge 3$, in which case the existence of non-trivial solutions in the half-space is open. The most interesting case seems to be $n=5$, when a non-trivial solution would give a finite-energy  boundary singularity for the steady-state equations. (As pointed out in [RF] and [St1], the 5d steady state problem is a good model for some aspects of the 3d time-dependent problem.)

\bigskip
{\bf 2. Regularity of very weak solutions.}

\medskip
By a
very weak solution of the steady-state Navier-Stokes system (\NSE)
in a domain $\Omega\subset R^n$ we mean a divergence-free vector
field $u=(u_1,\dots,u_n)\in L^2_{\hbox{\eightrm loc}}(\Omega)$ which
satisfies
$$
\int_{\Omega}
(u_i\Delta\varphi_i+u_iu_j{{\partial\varphi_i}\over{\partial
x_j}}) = 0
$$
for each smooth, compactly supported, divergence-free vector field
$\varphi=(\varphi_1,\dots,\varphi_n)$ in $\Omega$.

\medskip
 It is an open problem whether very weak solutions of
(\NSE) are regular. Standard regularity theory can be used to show
that very weak solutions are regular under the additional
requirement that $u\in L^n_{\hbox{\eightrm loc}}$ when $n\ge 3$. (In the case $n=2$ one can obtain regularity for $u\in L^{n+\ve}_{\hbox{\eightrm loc}},\,\,\ve>0$,  while the case $\ve=0, n=2$ appears to be open.) Equations
(\NSE) are usually considered with the assumption that $\nabla u\in
L^2_{\hbox{\eightrm loc}}$, in which case regularity follows for $n\le 4$ by a
standard bootstrapping argument. (The case $n=4$ is critical for $\nabla u\in L^2_{\rm loc}$.)  The assumption $\nabla
u\in L^2_{\hbox{\eightrm loc}}$ is of course very natural when
considering solutions describing real physical flows. However, one
can speculate that very weak solutions might arise from a blow-up
procedure of the usual weak solutions of the time-dependent
Navier-Stokes equations at a possible singularity (if a
singularity exists). The time-dependent 3-dimensional
Navier-Stokes equations are supercritical with respect to the
natural energy estimates, and in a blow-up procedure the
information about energy can be lost.

\medskip
A natural first step in understanding the regularity of the very
weak solutions above is to study the scale-invariant solutions in
$R^n$ which are smooth in $R^n\setminus{0}$. Theorem 3 in Section 6 settles this problem in $n\ge 4$.
For $n=3$ we can use Theorem 1: a calculation (which
can be found in [B], p.\ 209, and also in [T], and [CK]) shows
that, for $n=3$, Landau's solutions are not very weak solutions of (\NSE)
across the origin. Hence we have

\medskip
{\bf Corollary.} {\sl Let $n\ge 3$ and let  $u$ be a $(-1)$-homogeneous very weak
solution of the Navier-Stokes equations in $R^n$, which is smooth
away from the origin. Then $u\equiv 0$.}

\medskip
This result rules out only the simplest conceivable singularity of
a very weak solution. For example, the question if one can have a
non-trivial very weak solution smooth away from the origin and
satisfying $|u(x)|\le C|x|^{-1}$ in $R^n$ is not answered by
Theorem 1 and -- as far as I know -- remains open. (The results in [FR] can be used to
obtain some results for $n\ge 5$, under additional assumptions.)

\bigskip
{\bf 3. Landau solutions and behavior near $\infty$ of solutions in exterior domains in dimension $n=3$.}

\medskip
 Theorem 1 has some relevance for the problem of long-range
 behavior of solutions of the Navier-Stokes equations in three-dimensional exterior
 domains. (See for example [G] for an overview of this topic.)
 Let $f$ be a compactly supported vector field in $R^n$ and
 consider the equations
 $$
 \eqalign{-\Delta u +u\nabla u +\nabla p & = f \,,\cr
            \div u & = 0\,\cr}
             \quad\qquad\hbox{in $\,\,\,R^n$}\,\, , \eqno{\hbox{\rm (\Stn)}}
$$
together with a ``boundary condition'' at $\infty$, which might
take the form $u(x)\to 0$ at $\infty$ and $\int_{R^n} |\nabla u|^2
< \infty$, when $n\ge 3$. (See below for remarks concerning the case $n=2$, which is more difficult, essentially due to issues related to the Stokes paradox, see [Am].) The existence of such solutions was proved for $n=3$
in a classic paper by Leray ([Le]), but there are many open
questions about the behavior of these solutions for large $x$, see
[G]. The situation is more favorable in the case when $f$ is small,
 as in this case one can use perturbation techniques in spaces
 with prescribed decay for large $x$ to obtain a more detailed
 control of the solution.
 This technique was pioneered by R. Finn, see [F, G].
 However, the control is only in terms
 of the decay, it does not give the leading-order term, as the error
 term is of the same order of magnitude as the main term.

 Theorem 1 implies, roughly speaking,  the following:

\smallskip
{\bf Corollary.} {\sl In dimension $n=3$, if a solution of the above exterior problem is
asymptotically
$(-1)$-homogeneous, then the terms of order $|x|^{-1}$ must be
given by a Landau solution.}

\smallskip
To give this a more precise meaning, let us consider the scaled
functions $u_{\lambda}$ and $f_{\lambda}$ defined by $\ul(x)=
\lambda u(\lambda x)$ and $\fl(x)=\lambda^3f(\lambda x)$. The
functions $\ul$ and $\fl$ satisfy the same equations as $u$ and
$f$. Moreover the functions $\fl$ converge to a distribution
$\bar{ f}$, given by $\bar{f}(x)=b\delta(x)$, where
$b=\int_{R^3}f$ and $\delta$ is the Dirac function. {\it Assume}
now that $\ul$ converges to a limit $\bar{u}$ in, say,
$L^3_{\hbox{\eightrm loc}}(R^3\setminus\{0\})$ as
$\lambda\to\infty$. Our notion of ``asymptotically \hbox{$(-1)$-
homogeneous}" used above can be {\it defined} by requiring that
this is really the case. It is known that in the case of small
data this is true, see [NP].  The case of general large data remains open.  The limit functions $\bar{u}$ and $\bar{f}$
will again satisfy the same equations (in the sense of
distributions). Under our assumptions the function $\bar{u}$ is
smooth away from the origin, satisfies $\lambda\bar{u}(\lambda
x)=\bar{u}(x)$ for each $\lambda>0$, and, by Theorem 1, must
therefore be a Landau solution or vanish identically. (The
direction of the vector $b$ will be the axis of symmetry of the
solution.) For $b=0$ we will have $\bar u=0$, which means that,
under the above assumptions, the solution $u$ decays faster than
$|x|^{-1}$. After the first draft [Sv1] of this paper was written, a small data result similar to the above conclusions (for small data) was proved by a perturbation analysis
by A. Korolev and the author in [KS], without the use of Theorem 1 and [NP].

\medskip
The situation in dimension $n=2$ is quite different, and when $\int_{R^2} f\ne 0$, we do note expect existence of solutions to (\Stn) with $u(x)\to 0$ as $x\to\infty$, see for example [Am].
When $\int_{R^2} f=0$ we can write (under our assumptions) $f=\div F$ for a compactly supported smooth $F$. We will see in Section 5 that the  symmetry $u(x)\to \lambda u(\lambda x), F(x)\to\lambda^2 F(\lambda x)$ leads to a formal possibility of a limit given by \moh{} solution. However, this can happen only in exceptional cases, see the discussion in Section 5.

\medskip
The situation in dimensions $n\ge 4$ is simpler, at least at the formal level. In this case the leading term of the solution at $x\to \infty$ should be given by the linearized equation. For small data this can be indeed established by suitable perturbation arguments, whereas the large data situation has not been much studied, it seems.

\vbox{\bigskip
{\bf 4. Proof of Theorem 1.}

\medskip
Let $u$ be a $(-1)$-homogeneous
vector field in $R^3$, smooth away from the origin. Clearly $u$ is
determined by its restriction to the unit sphere $S^2\subset R^3$.
For $x\in\s2$ we decompose $u(x)$ as $u(x)=v(x)+f(x)e(x)$, where
$e(x)=x$ is the outer unit normal to $\s2$, and $v(x)$ is tangent
to $\s2$ at $x$, i.\ e.\ $v(x)\cdot n(x)=0$. We now write down the
Navier-Stokes equations for $u$ and as a system of PDEs on $\s2$.
If $u$ satisfies the Navier-Stokes equation in $R^3\setminus\{0\}$
in the very weak sense defined above, it is easy to see that there
exists a suitable pressure function $p$ in $R^3\setminus\{0\}$
which is $(-2)$-homogeneous and smooth away from the origin. The
function $p$ is also determined by its values on $S^2$, and the
system (\NSE) can be written down as a system of PDEs on $S^2$ for
$v, f$ and $p$. The differential operators in what follows will
all be differential operators on $S^2$, defined by the usual
conventions of Riemannian geometry. The differential forms on
$S^2$ will be identified with vector fields and vice-versa, as is
usual on Riemannian manifolds. The Hodge Laplacian $d\,d^*+d^*d$
on 1-forms will be denoted by $-\hl$. (The reason for writing it
as $-\hl$, with the minus sign, is to keep the equations on $\s2$
in a form which resembles the standard euclidean form of the
equations as much as possible.) The Navier-Stokes equations (\NSE)
for $u$ written in terms of $v,f$ and $p$ as equations on $S^2$
are as follows:}
$$ \eqalign{
   -\hl v +v\nabla v +\nabla (p-2f) & =0 \,,\cr
            -\Delta f + v\nabla f-f^2-|v|^2-2p & = 0\,,\cr
            \div v +f & = 0\,.} \eqno\hbox{\rm(\S)}
$$
A straightforward (although perhaps not the most illuminating) way
to derive these equations is to write the system (\NSE) in
spherical coordinates (see, for example, [B], p.\ 601) and check
that for $(-1)$-homogeneous vector fields it reduces to the system
(\S). We remark that the spherical coordinates version of (\NSE) in
the second edition of the book [LL]  (p.\ 49) contains a misprint
in the right-hand side of the first equation, where an incorrect
expression $\sin^2\theta$ appears instead of the correct
$\sin\theta$.
For the convenience of the reader we give another derivation of the
equations (\S) in Appendix 1.

\medskip
We will denote by $\omega$ the function on $S^2$  given by
$dv=\omega\Omega_0$, where $\Omega_0$ is the canonical volume form
of $S^2$. This corresponds to the formula $\omega=\curl v$ used in
$R^2$.

\medskip By taking $d$ of the first equation of the system (\S) we
obtain (see Appendix 2)
$$
-\Delta \omega + \div(v\omega) = 0. \eqno\hbox{\rm (\Eone)}
$$

\medskip
\vbox{ {\sl Lemma 1.} {\sl With the notation introduced above, we
have $\omega\equiv 0$.}

\smallskip
{\sl Proof.} Let $L$ be the differential operator defined by
$Lw=-\Delta w + \div(vw)$. The adjoint operator $L^*$ is given by
$L^*w=-\Delta w - v\nabla w$. The kernel of $L^*$ consists of
constant functions, as can be seen from the strong maximum
principle. The kernel of $L$ must therefore also be one
dimensional. Let us denote by $w_0$ a non-trivial function in the
kernel of $L$. If $w_0$ changed sign on $S^2$, we could find a
strictly positive smooth function $h$ on $S^2$ with
$\int_{\s2}w_0h=0$. But this would mean that the equation
$L^*w_1=h$ has a solution. However, the last equation cannot be
satisfied at  points where $w_1$ attains its minimum. From this we
see that the function $\omega$ cannot change sign. At the same
time, the definition of $\omega$ immediately implies that
$\int_{\s2}\omega=0$, and we see that $\omega$ must vanish.}

\smallskip
{\it Remark.} I assume the above argument is  known in one form or
another, but I was not able to find a good reference for it.

\medskip
Once we know that $dv=0$, the first equation of (\S) simplifies.
Indeed, when $dv=0$ we have $-\hl v=-\nabla\div v=\nabla f$, and
we also have $v\nabla v=\nabla |v|^2/2$. Using this, the first
equation of (\S) implies
$$
{1\over2} |v|^2+p-f=c\,,
$$
where $c$ is a constant. The second equation of (\S) now gives
$$
-\Delta f-2f+\div(fv)=2c\,. \eqno{\hbox{\rm (\Etwo)}}
$$
Integrating (\Etwo) over $S^2$ and using the third equation of (\S) we
see that $c=0$. Since $dv=0$ we can write $v=\nabla \varphi$ for a
suitable smooth function $\varphi$ on $\s2$. The equation (\Etwo),
together with the third equation of (\S) and the fact that $c=0$
now gives
$$
\Delta^2\varphi+2\Delta\varphi-\div(\Delta\varphi\nabla\varphi)=0\,.
\eqno{\hbox{ (\Ethree)}}
$$

Letting $w=2-\Delta\varphi$, the last equation can be re-written
as
$$
-\Delta w + \div (\nabla \varphi w) =0\,.
$$
The solutions of this equation are well-known: They are functions
of the form $c_1e^\varphi$, where $c_1$ is a constant. (An easy
way to verify this is for example the following: Write $w$ in the
form $c_1(x)e^{\varphi(x)}$. We get an equation for $c_1$  for
which the strong maximum principle implies that the solutions are
exactly $c_1(x)\equiv\hbox{\rm const.}$) Integrating $w$ over the
sphere we see that $c_1>0$. Hence we have
$$
-\Delta\varphi+2=c_1e^{\varphi}
$$
for a constant $c_1>0$. Changing $\varphi$ by a constant, if
necessary, we can assume $c_1=2$ without loss of generality, and
we end up with
$$
-\Delta\varphi+2=2e^{\varphi}\,.\eqno{\hbox{\rm (\Efour)}}
$$

\medskip
The interpretation of equation (E4) is well-known (see, for
example, [CY]): Let $\bar g$ be the canonical metric on $\s2$ and
let $g$ be the metric on $\s2$ defined by $g=e^\varphi\bar g$.
Equation (\Efour) says exactly that the Gauss curvature of the metric
$g$ is $1$, i.\  e.\ the metric $g$ is isometric to the metric
$\bar g$. In other words, we have $g=h^*\bar g$ (pullback of $\bar
g$ by $h$) for a suitable diffeomorphism $h$ of $\s2$. From the
definitions we also see that $h$ has to be conformal or
anti-conformal. Anti-conformal maps can be obtained from conformal
maps by a composition with an isometry, and hence we can only
consider the case when $h$ is conformal. For a given conformal
$h$, the function $\varphi$ is given by
$$
\varphi(x)=\log|h'(x)|^2\,, \eqno{\hbox{\rm (\Efive)}}
$$
where $h'(x)$ denotes the (complex) derivative of $h$ at $x$. It
is well-known (see e.\ g.\ [DFN]) that all conformal
diffeomorphisms of $\s2$ can be produced as follows. Let $P\colon
S^2\to\bf C$ be the standard stereographic projection, and let
$M_{\lambda}\colon\bf C\to\bf C$ be defined by $z\to \lambda z$.
Let $h_\lambda=P^{-1}\circ M_{\lambda}\circ P$. Then any conformal
diffeomorphism of $\s2$ can be produced by composing a suitable
$h_\lambda$ (with $\lambda>0$) with isometries of $S^2$. If
$\varphi$ is given by (\Efive) and we compose $h$ with and isometry,
then the function $\varphi$ either does not change or only changes
by being shifted by the isometry. Therefore in a suitable
coordinate frame all solutions $\varphi$ of (\Efour) look like the
solutions generated by the special $h_{\lambda}$ above. We now
consider the standard spherical coordinates $(\theta,\psi)$ on
$S^2$, given by
$$
\eqalign{
   x_1  = & \sin\theta\cos\psi \,,\cr
   x_2  = & \sin\theta\sin\psi\,,\cr
   x_3 = & \cos\theta\,.} \eqno\hbox{\rm(\SC)}
$$
We will use the usual notation $e_{\theta}={{\partial
x}\over{\partial \theta}}$ for the tangent vector field on $S^2$
corresponding to ${\partial}\over{\partial\theta}$. Letting
$\lambda=e^{-\kappa}$, calculating the maps $h_{\lambda}$ above in
these coordinates, and using the formula (\Efive), we obtain
$$
\varphi(x)=-2\log\,(\cosh\kappa-\sinh\kappa\,\cos\theta)\,.
\eqno{\hbox{\rm (\Esix)}}
$$
This gives
$$
\eqalign{
   v  = & \,\,{{\partial \varphi}\over{\partial \theta}}\,e_{\theta}={{-2\sin\theta}\over{\coth\kappa-\cos\theta}}\,e_\theta \,,\cr
   f = & -\Delta\varphi=2e^\varphi-2={2\over{(\cosh\kappa -
   \sinh\kappa\,\cos\theta)^2}}-2\,,} \eqno{\hbox{\rm (\Eseven)}}
$$
which agrees with the formulae in [B], p.\ 207 if we set $\coth
\kappa = 1+c$ and with the formulae in [LL], p.\ 82, if we set
$\coth \kappa=A$. The proof of Theorem 1 is finished.

\medskip
{\it Remarks:}

\smallskip
1. As we already mentioned in Section 1, the Landau solutions
(given by (\Eseven)) do not satisfy the Navier-Stokes equations (\NSE)
across the origin. A calculation in [B], p.\ 209, shows that for
Landau's solutions we have, in distributions,
$$
-\Delta u + \div (u\otimes u)+\nabla p = b\delta\,,
$$
where $\delta$ is the Dirac function and $b=b(\kappa)$ is a
non-zero vector in $R^3$ depending in a non-trivial way on the
parameter $\kappa$ which parametrizes the solutions in the above
coordinate frame. The exact formula for $b$ can be found in [B],
p.\ 209, and was also calculated in [CK].

\smallskip
2. If $h\colon S^2\to\s2$ is a non-trivial holomorphic map (which
is not necessarily a diffeomorphisms) the formula (\Efive) gives a
function $\varphi$ which is regular away from a finite set
$a_1,\dots,a_m\in S^2$ where $h'$ vanishes. The function $\varphi$
will generate a $(-1)$-homogeneous solution of the Navier-Stokes
equations in the region $R^3\setminus(\cup_{j=1}^{j=m}R_+\cdot
a_j)$, where $R_+=[0,\infty)$. However, the vector field will not
be locally square integrable in $R^3$, except for the case of
Landau's solutions, when $h'$ does not vanish at any point.

\bigskip
{\bf 5. (-1) - homogeneous solutions in dimension $n=2$.}

\medskip
In dimension $n=2$ the equations derived in Appendix 1
reduce to the circle $S^1$. We denote by $\theta$ the
the natural angle variable on the circle. The unknown
function are $v=v(\theta)$, the  component of the velocity
tangent to the circle, $f=f(\theta)$, the component of the
velocity normal to the circle, and the pressure
$p=p(\theta)$ on the circle. We will use the notation
$f'={d\over{d\theta}}f$. The equations are

$$ \eqalign{
   (p-2f)' & =0 \,,\cr
            - f'' + v f'-f^2-|v|^2-2p & = 0\,,\cr
            v'  & = 0\,.} \eqno\hbox{\rm(\Stwo)}
$$
This means that $v$ has to be constant and $p=2f+\hbox{\rm const.}$
Since we are looking for solutions on the whole circle, corresponding
to the periodic solutions in $\theta$ and the term $vf'$ can be interpreted
as ``damping", we see that $v$ or $f'$ must vanish identically.
The solutions corresponding to a non-zero $v$ are therefore the solutions
for which all unknown functions $f,v,p$ are constant and the constants satisfy
$f^2+v^2+2p=0$.

In the case of $v=0$ we obtain a single equation for $f$
$$
f''=-4f-f^2+b\,,
$$
where $b$ is any constant. This is the equation of motion of a particle
in the potential $V(f)=1/3 f^3 + 2f^2-bf$, and we are interested in
its $2\pi$-periodic solutions. From this interpretation and the form
of the potential $V$ it is clear than one has many of such solutions.
(The key point is that for large $b$ the potential has a local
minimum, and the solutions of the linearization of our equation
around this equilibrium oscillate at high frequency. By changing the
amplitude of the oscillations we can change the period and adjust it
so that the solution is periodic with the smallest period $2\pi/m$
for a positive integer $m$. Together with the freedom to change
$b$, this gives countably many 1-parameter families of solutions.
We refer the reader to [Ha] for the details.)

We will be interested in the solutions which satisfy the additional
requirement that
$$
\int_{S^1} f = 0\,. \eqno\hbox{\rm (\C)}
$$
In dimensions $n\ge 3$ this condition is satisfied automatically
due to the equation $(n-2)f=-\div v$. (This is also reflected by the fact that
in dimensions $n\ge 3$ any vector field $u$ in $R^n$ which is
div-free in $R^n\{0\}$ and bounded by $c/|x|$ is also div-free
in $R^n$ in the sense of distribution. In dimension $n=2$ this is
no longer the case.)
Condition (\C) comes up naturally in the context of
the long-distance behavior of steady-state solutions in the following way.
For a matrix field $F=F_{ij}(x)$ in $R^n$ we denote by
$\div F$ the vector field ${\partial\over{\partial x_j}}F_{ij}$.
We consider the steady Navier-Stoked equations with the right-hand side
in the divergence form:
$$ \eqalign{-\Delta u +u\nabla u +\nabla p & =\div F \,,\cr
            \div u & = 0\,.\cr} \eqno\hbox{\rm(\NSEdiv)} $$
In dimension $n=2$ the scaling symmetry works for this equation
in a way similar to the 3d case with $\div F$ replaced by $f$,
which we dealt with in Section 3: for $\lambda>0$ the quantities
$$ \eqalign
{u_\lambda(x) & =\lambda u(\lambda x)\,,\cr
            p_\lambda(x)& =\lambda^2 p(\lambda x)\,, \quad\hbox{\rm and}\cr
            F_\lambda(x) & =\lambda^2 F(\lambda x)} \eqno\hbox{\rm (\SCtwod)}
            $$
satisfy again equation (\NSEdiv).

Assume now that $\lim_{\lambda\to\infty} u_\lambda=\bar u$ exists.
Then $\bar u$ is $(-1)$-homogeneous. It obviously satisfies
$\int_{S^1} u\cdot\nu=0$, where $\nu$ is the outer unit normal
to $S^1$. In the variables $(f,v)$ above this means that
$\int_{S^1} f = 0$.
Clearly $\lim_{\lambda\to\infty} F_\lambda = \bar F = M\delta$ where
$M=\int_{R^2} F$ is a $2\times 2$ matrix and $\delta$ is the Dirac function, and formally
one has
$$
-\Delta\bar u + \bar u\nabla \bar u + \nabla\bar p = \div\bar F
$$
for a suitably defined $\bar p$. For the linear Stokes problem
(obtained by dropping the term $u\nabla u$ from the equations)
the above procedure works well and the field $\bar u$ gives the
leading terms asymptotics at $\infty$ for the solution.

The solutions of
$$
-\Delta \bar U  + \nabla \bar P = \div (M\delta)
$$
are given by
$$
\bar U_i(x)=M_{jk}G_{ijk}(x)\, ,
$$
where
$$
G_{ijk}(x)={1\over{4\pi}}{\partial\over{\partial x_k}}\left(\delta_{ij}\log{1\over |x|}+{{x_ix_j}\over{|x|^2}}\right) \,,
\eqno\hbox{(\Gtwo)}
$$
and these solutions give the leading-order behavior of the solutions of
$$
-\Delta  U  +  \nabla P = \div F \eqno{\hbox{(\St)}}
$$
as $x\to\infty$.
One way to calculate the Green function (\Gtwo) is to solve the  linearization
of (\Stwo), which can be easily done explicitly. For example, the
vector field $G_{i11}$ corresponds to $v=0, \,\,f(\theta)={1\over{4\pi}} \cos(2\theta)$.

In dimension $n=3$, with $\div F$ replaced by $f$ this procedure works also
at the non-linear level, at least for small data, as we have seen in Section 3.

Can this also work for the non-linear problem in dimension $n=2$?
One difficulty is
that for $(-1)$-homogeneous functions the term $\bar u\nabla \bar u$ no longer has an easy
distributional interpretation in the open sets containing $x=0$.
Even if we write is as $\div u\otimes u$, it is still not transparently
well-defined as a distribution.
We can side-step this issue by considering the equations only in $R^2\setminus\{0\}$.
The functions $\bar u$ in $R^2\setminus\{0\}$ which are results of the
above ``blow-up procedure" will still satisfy the equation $\div \bar u = 0$
across the origin, which translates to $\int_{S^1} f = 0$.
We see that it important  to characterize the solutions of  (\Stwo) which satisfy
the zero flux condition (\C).
These solutions are characterized in the following theorem.

\medskip
\vbox{
{\bf Theorem 2.} {\sl The solutions of (\Stwo) satisfying the zero flux conditions
(\C) are of the following form:

Either
$$
f=0,\quad v=\hbox{\rm const.}, \quad \hbox{and $p=-|v|^2/2$,}
$$
or
$$
f(\theta)= \tilde f_k(\theta-\theta_0), \quad v=0, \quad\hbox{and $p=-f/2+c_k$},\quad k=3,4,\dots\,,
$$
where for each $k=3,4,\dots$, the function $\tilde f_k$ is a  non-trivial periodic function of $\theta$ with minimal period $2\pi/k$, the constant $c_k$ is given by $c_k={1/2}\int_{S^1}|\tilde f_k|^2$, and $\theta_0$ can be chosen
 arbitrarily. The functions $\tilde f_k$ can be expressed in terms of
the classical elliptic functions. The amplitude of oscillations of $\tilde f_k$ is of order $k^2$. }
}

\medskip
Before going to the proof of the theorem, let us point out an interesting conclusion one can
make from it. Let us consider the equation (\NSEdiv) in dimension $n=2$ with a smooth
compactly supported $F$ satisfying
$\int_{R^2} F_{11}\ne 0$, and $F_{12}=F_{21}=F_{22}=0$. (One can say that the force $\div F$
is approximately a dipole in the $x_1$-direction.)
The solution of the linear Stokes system (\St) is given by
$$
U_i=G_{i11}*F_{11}
$$
and its asymptotics as $x\to \infty$ is given by (a multiple of) $G_{i11}$, modulo terms of order $1/|x|^2$.
The field $G_{i11}$ (the solution corresponding to an exact unit ``dipole force" in the $x_1$-direction)
corresponds to the solution of the linearization of the system (\Stwo) with $f=\cos(2\theta)$.

\smallskip
One can now ask if in the situation when $F$ is small, one has a solutions of the full Navier-
Stokes equation with a similar structure. Theorem 2 shows that, somewhat surprisingly, this
is not the case: the system (\Stwo) does not have any  solution which would be close to the
solution $f=\cos(2\theta),\, v=0$ of the linearized system. Therefore the linear solution
$G_{i11}*F_{11}$ cannot be ``deformed" into the solution of the full non-linear system
which would still be asymptotically $(-1)$-homogeneous as $x\to\infty$, no matter
how small $F$ is, as long as $\int_{R^2} F_{11}\ne 0$. In particular, one cannot
obtain solutions of (\NSEdiv) for small $F$ by perturbation techniques in the spaces
of functions with decay $O(1/|x|)$ as $x\to\infty$. The failure of the usual perturbation
series to converge in the spaces with decay $O(1/|x|)$ can be analyzed in some detail
and is interesting by itself. A noteworthy feature of the situation is
that the failure does not occur at the level of the ``second iterant" (with the first iterant
being the linear solution), but only at the level of the third iterant.
Some solutions of (\NSEdiv) can be constructed by Leray's method based on solving
the problem in large balls $B_R$ by using energy estimates together with some topological
arguments (e.\ g. degree theory), and then letting $R\to\infty$. However, the precise behavior
for large $x$ of the solutions obtained in this way seems to be open.

\medskip
{\bf Proof of Theorem 2}.
As we have already mentioned, the problem without the zero-flux condition (\C) has been investigated
in some detail in 1917 by G. Hamel, [Ha]. For the proof we will change our notation and
instead of $f=f(\theta)$ we will write $u=u(\theta)$ for the radial component. It is clear
that the only non-trivial part of the proof is the investigation of the solutions with $v=0$.
This reduces our task to problem of find all non-trivial $2\pi$-periodic solutions of
$$
u''=-4u-u^2+b  \eqno{\hbox{(\ODEone)}}
$$
with $\int_0^{2\pi} u(\theta)\,d\theta = 0$, where $b$ is an arbitrary real parameter.
As above, we interpret the solutions as motions of a particle of unit mass in the
potential $V(u)=u^3/3+2u-bu$. Therefore we have the usual energy conservation
$$
(u')^2=2E-2V(u) \eqno{\hbox{(\E)}}\,.
$$
This is a classical equation defining the elliptic functions (see, for example, [Ch]).
Following [Ha], we note that the  relevant situation for us occures exactly when
the polynomial  $2E-2V(u)$ has three real roots $e_1\ge e_2\ge e_3$ satisfying
$$
e_1+e_2+e_3=-6\,. \eqno{\hbox{(\rc)}}
$$
(Instead of choosing $b$ and $E$ we choose the roots
$e_i$ satisfying (\rc).)
One can therefore write
$$
u'=\pm\sqrt{{2\over3}(e_1-u)(u-e_2)(u-e_3)}\,,
$$
and our task is to investigate for which choices of the roots we have
$$
T=\int_{e_2}^{e_1} {du\over{\sqrt{(e_1-u)(u-e_2)(u-e_3)}}}=\sqrt{2\over 3}{\pi\over k}\quad \hbox
{for some $k=1,2,\dots$}  \eqno{\hbox{(\P)}}
$$
(which says that $u$ is $2\pi$-periodic) together with
$$
I=\int_{e_2}^{e_1}{u\,du\over{{\sqrt{(e_1-u)(u-e_2)(u-e_3)}}}}=0\,, \eqno{\hbox{(\Cprime)}}
$$
(which is just another way of stating condition (\C).)

Following [Ha], we use the classical change of variables in these elliptic integrals:
$$
u=e_2+(e_1-e_2)\sin^2\varphi\,,
$$
and we also set
$$
\kappa={{e_1-e_2}\over{e_2-e_3}}, \quad \delta = e_2-e_3\,.
$$
and
$$
F(\kappa)=\int_0^{\pi/2}{d\varphi\over{\sqrt{1+\kappa\sin^2\varphi}}},\quad E(\kappa)=\int_0^{\pi/2}\sqrt{1+\kappa\sin^2\varphi}\,d\varphi\,.
$$

This gives
$$
T={2\over{\sqrt{\delta}}}F(\kappa) \quad \hbox{and $\,\,$ ${\,\,\sqrt{\delta}\over 2}I=-2-{\delta(2+\kappa)\over{3}}F(\kappa)+ \delta E(\kappa)$}\,.
$$
The functions $F(\kappa)$ and $E(\kappa)$ are variants of the classical complete elliptic
integrals of the first and second kind, respectively. (In the classical definition one
replaces $\kappa$ by $-m=-k^2$ and $F$ is denoted by $K$.)
It is easy to check that the triples of roots with $e_1>e_2>e_3$ and $e_1+e_2+e_3=-6$ are in one-to-one
correspondence with the pairs $\kappa>0,\,\delta>0$.

\smallskip
We use the condition $I=0$ to obtain
$$
{2F\over\delta}=E-{2+\kappa\over 3}F\,.
$$

We we see
 that $I=0$ can only be satisfied when $E(\kappa)-{1\over3}(2+\kappa)F(\kappa)\ge 0$ and
in that case we have
$$
T^2={4F^2\over\delta}=2F(\kappa)\left(E(\kappa)-{2+\kappa\over 3}F(\kappa)\right)\,.
$$
Let us denote by $H(\kappa)$ the function on the right-hand side.
We are interested in the non-negative solutions of the equation
$$
H(\kappa)={{2\pi^2}\over{3k^2}}\,, \quad k=1,2,\dots
$$
We note that $H(0)={\pi^2\over 6}={{2\pi^2}\over{3\cdot2^2}}$. The solution $\kappa=0$ corresponds to
the roots $e_1=e_2=0$ and $e_3=-6$ and the ``infinitesimal oscillations" of $u$ around $u=0$, which
is exactly the solution of the linearized equation. Its period is $\pi$, as expected.
From the definitions of $F$ and $E$ it is also easy to see that $H(\kappa)$ becomes
negative for sufficiently large $\kappa>0$. Hence the proof of the existence of $\tilde f_k$
and their uniqueness (modulo the shift by $\theta_0$)
will be finished if we show
that the derivative $H'(\kappa)={dH(\kappa)\over d\kappa}$ is strictly negative for $\kappa>0$.

\smallskip
We will need the classical formulae for the derivatives of $E,F$
$$
\eqalign{
E' = {d E\over d\kappa} & ={1\over{2\kappa}}(E-F) \,,\cr
            F'={dF\over d\kappa} & = {1\over{2\kappa}} ({E\over{1+\kappa}}-F)\,
            }
            \eqno\hbox{\rm(\der)}
$$
together with the inequality
$$
1< {E\over F}<1+{\kappa\over 2}, \quad \kappa>0, \eqno{\hbox{(\efi)}}.
$$
see Appendix 3.

We calculate
$$
3\kappa(1+\kappa) H'=3E^2-2(2+\kappa)EF+(1+\kappa)F^2=F^2(3x^2-2(2+\kappa)x+(1+\kappa)), \quad x={E\over F}\,.
$$
It is not hard to see that
$$
3x^2-2(2+\kappa)x+(1+\kappa)<0, \quad \hbox{when $\kappa>0$ and $1<x<1+{\kappa\over 2}$},
$$
which shows that $H'(\kappa)<0$ for $\kappa>0$.

 The amplitude of the oscillation of $\tilde f_k$ is
$$
e_1-e_2=\kappa\delta={4\kappa F^2\over H}={6\kappa F^2\over\pi^2}k^2\,,
$$
which proves the statement about the amplitude of $\tilde f_k$, as for $k\to\infty$ the corresponding
values of $\kappa$ converge to the positive root of the equation $H(\kappa)=0$.
This finishes the proof of Theorem 2.

\bigskip
{\bf 6. Higher Dimensions.}

\medskip
In this section we show that the system (\Sn) does not have (smooth) solutions in dimensions $n\ge 4$.
As we already indicated in Section 2, this  result is related to the regularity theory of
the steady-state equations. A $(-1)$-homogeneous solution in dimension $n=5$ would provide the
simplest example of a singular solution with locally finite energy $\int_{B_R}|\nabla u(x)|^2\,dx$.
Dimension $n=5$ is the lowest dimension for which the steady Navier-Stokes is super-critical
with respect to the energy $ \int_{B_R}|\nabla u(x)|^2\,dx$, in the sense that the classical
boot-strapping argument cannot be used to prove regularity. The regularity theory for this case
has been studied by Frehse and R\accent23 u\v zi\v cka, see for example [FR],
and by Struwe, see [St1]. The key point of these works is to use the special properties
of the quantity $|u|^2/2+p$. The quantity will also play an important role in the
proof Theorem 3 below which is the main result of this section. The theorem follows
from the results of Frehse and R\accent23 u\v zi\v cka, and was also proved
by Struwe [St2],  and T.-P. Tsai and the author, see [T]. For the convenience of the
reader we reproduce below the proof by T.-P. Tsai and the author.

\bigskip
{\bf Theorem 3.} {\sl When $n\ge 4$, the system (\Sn) has no non-trivial solutions.}

\bigskip
{\bf Proof.} The key point in the proof is to use a well-known non-trivial identity which
is satisfied by the ``Bernoulli quantity"  $H=|u|^2/2+p$ for any steady-state Navier-Stokes
solutions. Denoting by $\omega$ the anti-symmetric part of $\nabla u$, we have
$$
-\Delta H + u\nabla H = -2|\omega|^2\,. \eqno{\hbox{(\B)}}
$$
This identity plays a very important role in the regularity theory for higher-dimensional
steady-state Navier-Stokes. For a $(-1)$-homogeneous solution we  will denote, with a slight abuse of notation,
by $H, |\omega|^2,$ and $p$ also the restriction of these quantities (originally defined in $R^n\setminus\{0\}$)
to the sphere $S^{n-1}$. We recall that we write the restriction of the vector field $u$ to $S^{n-1}$
as $u=v+fe$, where $v$ is tangential to the sphere and $e$ is the normal to the sphere.
For the proof of Theorem 3 it is enough to replace the first equation of (\Sn) by
the equation (\B) expressed in terms of the variables on $S^{n-1}$. This system is
$$ \eqalign{
   -\Delta H +(2n-8)H +v\nabla H -2fH & =-2|\omega|^2 \,,\cr
            -\Delta f + v\nabla f  & = 2H\,,\cr
            \div v +(n-2)f & = 0\,,} \eqno\hbox{\rm(\SnB)}
$$
where all the differential operators are now taken on $S^{n-1}$.
When $n=4$, we can integrate the first equation over the sphere.
Integrating by parts and using the third equation we see that
the integral of the left-hand side vanishes, and hence $\omega$
must vanish identically.
When $n>4$, we let $H_+$ be the positive part of $H$ and $\alpha=(n-4)/2$.
We multiply the first equation by $H^\alpha_+$ and integrate by parts to obtain
(with the use of the third equation)
$$
\int_{S^{n-1}}\left(\alpha|\nabla H|^2H^{\alpha-1}_+  +(2n-8)H^{1+\alpha}_+ +|\omega|^2H^{\alpha}_+\right)\,dy = 0\,.
$$
This shows that $H_+$ has to vanish, which means that $H$ is non-negative.
Now the second equation of (\SnB) together with the strong maximum principle imply
$f$ must be constant, and therefore also $H$ vanishes. Going back to the fist
equation we see that we obtain again that $\omega$ must vanish identically.
We now look again at the $(-1)$-homogeneous field $u$ defined in $R^n\setminus\{0\}$.
Since $\omega=0$ and $\div u=0$, we see that $u$ is harmonic in $R^n\setminus\{0\}$,
and since $n\ge 4$, the $(-1)$-homogeneous singularity at $x=0$ is removable.
Hence $u$ vanishes identically.

\medskip
{\bf 7. Open problems.}

\medskip
 An interesting problem is to try to
repeat, the above analysis  when $\rmo n$ is replaced by the half-space
$R^n_+=\{x\in R^n, \, x_n>0\}$ and the boundary condition $u=0$ is
imposed on $\partial R^n_+\setminus\{0\}$. The case $n=2$ is amenable to an ODE analysis, along the lines of Section 5, see also [Ha].
When $n\ge 3$, the problem becomes more difficult, and the following question seems to be open.

\smallskip
{\sl For $n\ge 3$, does Theorem 3 remain true for in $R^n_+\setminus\{0\}$, with the boundary condition $u=0$ at $\partial R^n_+\setminus\{0\}$ ?}

\smallskip
If an analogue of Theorem 3 would fail in dimension $n\ge 5$ and a non-trivial solution existed, one would have a genuine example of a boundary singularity for steady-state solutions with (locally) finite energy $\int_{B_r\cap R^n_+} |\nabla u|^2$ (in the corresponding dimension).

\smallskip
In dimension $n=3$, a relatively simple calculation
shows that there are no non-trivial axi-symmetric
$(-1)$-homogeneous solutions in that case. However, it
is not clear whether this conclusion is still true without assuming
the rotational symmetry. We refer the reader to the very
interesting paper [Se], where a related situation is studied in a
different context.

\medskip
Another interesting question is the following:

\smallskip
{\sl Among smooth
vector fields in $R^3\setminus\{0\}$ satisfying $|u(x)|\le
C|x|^{-1}$ for some $C>0$, are the Landau solutions the only ones
which satisfy the Navier-Stokes equations (\NSE) in
$R^3\setminus\{0\}$?}

\smallskip
Such questions are relevant for the problem
of asymptotic behavior of steady-state solutions in exterior
domains mentioned in Section 3. A first natural step in addressing
this question is to look at possible infinitesimal deformations of
Landau solutions in the above class. This leads to linear
equations which can be reduced to ODEs by classical methods of
separation of variables, due to the symmetries of Landau's
solutions. Based on numerical experiments with these ODEs, the
author conjectures that the Landau solutions are rigid with
respect to infinitesimal deformations, i.\ e.\ it seems that there
are no new solutions bifurcating from Landau's solutions.

\bigskip
\bigskip
{\bf Appendix 1.}

\medskip
In this section collect some formulae which can be used
for an alternative derivation of equations (\S) and (\Eone).
As we mentioned in Section 4, (\S) and (\Eone) can be checked by
straightforward  but tedious calculations in polar coordinates.
However, it seems to be useful to have a more illuminating
derivation.

\medskip

Let us consider a $(-1)$-homogeneous vector field $u$ in $R^n$
which is smooth away from the origin. We will write the
coordinates in $R^n$ as $x=(x_1,\dots,x_n)$, and denote
$r=|x|$ the distance to the origin. We can write
$$
x=ry\,,
$$
 with $y\in S^{n-1}$, where $S^{n-1}\subset R^n$ is the standard
 unit sphere. For $y\in S^{n-1}$ we let $e(y)=y\in R^n$ be the outward unit normal.
 The vector field $u$ can be written as
 $$
 u(x)={1\over r}\left(v(y)+f(y)e(y)\right)\,,
 $$
 where $v$ is a vector field tangent to the sphere and $f$ is a function
 on the sphere.

 \medskip
 We would like to express the Navier-Stokes equations (\NSE) for $u$ in terms
 of intrinsic equations on $\sm$ for the field $v$ the function $f$ and the pressure.
 It is easy to see that the pressure (which is only given up to a constant) can be
 chosen so that
 $$
 p(x)={1\over {r^2}}p(y)\,.
 $$
 The function $p(y)$ can then be considered as function on $\sm$.

 \medskip

 If we dealt with the Euler equations
 $$ \eqalign{ u\nabla u +\nabla p & =0 \,,\cr
            \div u & = 0\,.\cr} %\eqno\hbox{\rm(NSE)}
            $$
rather than the Navier-Stokes, the derivation would be straightforward:
we would get

$$ \eqalign{
   v\nabla v +\nabla p & =0 \,,\cr
             v\nabla f-f^2-|v|^2-2p & = 0\,,\cr
            \div v +f & = 0\,,} %\eqno\hbox{\rm(\S)}
$$
where all the differential operators are the intrinsic operators on $\sm$. For example,
$v\nabla v$ is the covariant derivative of $v$ in the direction of $v$.
This calculation follows directly from the definition of the covariant derivative
in terms of the ``usual derivative" and the orthogonal projection on the
tangent space, and it is left to the reader as an easy exercise.

\medskip

For Navier-Stokes we must include the Laplacian $\Delta u$ and expressing
this term in suitable intrinsic operators on $\sm$ is more subtle, although
such calculation are routine in Differential Geometry.

\medskip
We will consider the following operators:

\smallskip
$\difsm$ is the standard differentiation (of $R^k$-valued functions, $k=1,2,\dots$) on $\sm$,

\smallskip
$\difc$ is the covariant differentiation of the vector fields (or one-forms) on $\sm$,

\smallskip
$\lrn$ is the standard Laplacian on $R^n$, corresponding to the quadratic form $\int_{R^n} {1\over 2} |\nabla X|^2$

\smallskip
$\lsm$ is the standard Laplacian (on $R^k$-valued functions, $k=1,2,\dots$) on $\sm$, corresonding
to the quadratic form $\int_{\sm}\pul |\difsm X|^2$

\smallskip
$\cl$ is the covariant Laplacian (also called ``rough Laplacian") on vector fields or one-forms on $\sm$,
corresponding to the quadratic form $\int_{\sm} \pul |\difc X|^2$

\smallskip
$\hl$ is the Hodge Laplacian on vector fields or one-forms on $\sm$, corresponding to the quadratic form
$\int_{\sm} \pul (|d X|^2 + |d^* X|^2)$, where $d$ is the exterior differentiation and $d^*$ its adjoint
(essentially the operator div).

\smallskip
$\ric$ is the Ricci curvature tensor on $\sm$. Recall that ${\ric}=\{R_{ij}\}_{i,j=1}^{(n-1)}$, and
$R_{ij}=(n-2)g_{ij}$, where $g_{ij}$ denotes the metric.

\medskip
We also recall the formula
$$
\int_{\sm} \left(|\difc u|^2+{\ric}(u,u)\right) =\int_{\sm} (|du|^2 + |d^*u|^2)\,\,, \eqno\hbox{(\Aone)}
$$
which follows by integration by parts. This formula implies the identity
$-\cl+\ric = -\hl$. Taking into account that we are on $\sm$, we can write
$-\cl+(n-2)=-\hl$.

\medskip
We also recall that
$$
\Delta = {{\partial ^2}\over{\partial r^2}}+{{(n-1)\partial}\over{r\partial r}}+{1\over{r^2}}\lsm\,.\eqno\hbox{(\Atwo)}
$$

Therefore, returning to our $-1$ -homogeneous field $u={1\over r}(v(y)+f(y)e(y))$,  we have
$$
\Delta u = {1\over r^3}\left((3-n)(v(y)+f(y)e(y))+\lsm(v(y)+f(y)e(y))\right)\,\,.
$$
Here the  Laplacian on  the right-hand side is the usual sphere Laplacian of the
$R^n$-valued function $v(y)$ on $\sm$, i.\ e.\ we calculate it
``component by component". We need to decompose this expression into
the tangent part and the normal part, and write each part in terms
of intrinsic operators on on vector-fields/one forms and functions
on the sphere.  An easy way to do this is to use the corresponding
quadratic forms. Let $X(y)=v(y)+f(y)e(y)$ and let us consider
the quadratic form
corresponding the $\lsm X$, which is
$$
\int_{\sm}\pul|\difsm X|^2
$$

For a fixed vector $b$ tangent to the sphere we have
$$
\difsm_b X = \difc_b v - II(b,v)e+(\difc_b f)e+fb\,,
$$
where $II(b,v)$ denotes the second fundamental form, which in our
case is simply the scalar product $(v,b)$. Evaluating $|\difsm_b X|^2$
and summing over orthonormal vectors $b$, we obtain
$$
|\difsm X|^2=|\difc v|^2 + 2 f\, \div v + (n-1)|f|^2 + |v|^2 -2 v\difsm f + |\difsm f|^2\,.
$$
Integrating this identity over the sphere and using (\Aone) we see that
$$
\int_{\sm}|\difsm X|^2=\int_{\sm}\left(|dv|^2+|d^*v|^2+ (3-n)|v|^2+|\difsm f|^2+(n-1)|f|^2+4f\,\div v\right)\,.
$$
Taking variations of the form $\delta X=\varphi(y)+\eta(y)e(y)$ with $\varphi(y)$ tangent
to the sphere, we see that the tangential part of $-\lsm X$ is
$$
[-\lsm(v+fe)]_{\hbox{\fiverm tangential}} = -\hl v+ (n-3)v-2\difsm f
$$
and the normal part is
$$
[-\lsm(v+fe)]_{\hbox{\fiverm normal}} = -\lsm f + (n-1)f + 2\div v \,.
$$
We recall that the continuity equation $\div u = 0$ implies $\div v = -(n-2)f $ and hence we can write
$$
[-\lsm(v+fe)]_{\hbox{\fiverm normal}} = -\lsm f + (3-n)f.
$$
Using (\Atwo) together with
$$
[{{\partial ^2}\over{\partial r^2}}+(n-1){{\partial}\over{r\partial r}}]\,\, {1 \over r}= {(3-n)\over {r^3}}\,,
$$
we arrive at
$$
\eqalign{
[-\Delta u]_{\hbox{\fiverm tangential}} & = {1\over{r^3}}(-\hl v - 2\difsm f)\,,\cr
[-\Delta u]_{\hbox{\fiverm normal}}  & = {1\over {r^3}} (-\lsm f)\,.
}
$$

Putting this together with the Euler part above, and dropping the indices S and C in $\difsm, \difc, \lsm$ since
all equations are now intrinsic on the sphere and there is no danger of confusion,
we see that the Navier-Stokes for $u$ becomes
the following system on $\sm$:

$$ \eqalign{
   -\hl v +v\nabla v +\nabla (p-2f) & =0 \,,\cr
            -\Delta f + v\nabla f-f^2-|v|^2-2p & = 0\,,\cr
            \div v +(n-2)f & = 0\,.} \eqno\hbox{\rm(\Sn)}
$$

\bigskip
\vbox{{\bf Appendix 2.}

\medskip
We consider equation (\Eone), which was obtained in Section 4
from the first equation of (\S) by applying $d$.
We recall that $\Omega_0$ denotes the volume form on $S^2$,
and that $\omega=\curl v$ is defined by $dv=\omega\Omega_0$.
Both $v$ and $v\nabla v$ can be considered  as vector
fields or one-forms, depending on the context.
For a vector field $X$ we will denote by  $L_X$ the Lie derivative
along $X$.}

The formula which we would like to prove is
$$
d(v\nabla v)= \Omega_0\,\,\div(v\omega)
$$
or, equivalently,
$$
\curl (v\nabla v)=\div(v\omega).
$$
This equation explains the special behavior of vorticity in two dimension.
It can verified by mechanical calculation. However, we prefer a more
geometric derivation, which avoids the calculations and gives a better
explanation of this identity, even in the flat case.
We will use the traditional notation $v^i_{,j}$ and $v_{i,j}$ for
covariant differentiation. We have
$$
(v\nabla v)_i=v^j v_{i,j}=v^jv_{i,j}+v^jv_{j,i}-v^jv_{j,i}=(L_v v)_i-{1\over 2}(v^jv_{j,i}+v_jv^j_{,i})=(L_v v)_i-{1\over2}(v_jv^j)_{,i}\,.
$$
As it is hopefully clear from the context, in the expression $L_v v$ the first $v$ is considered as a vector field, whereas the second $v$ is considered as a one-form.
We see that
$$
d(v\nabla v)=d(L_v v)=L_v(dv)=L_v(\omega\Omega_0)%(L_v\omega)\Omega_0+\omega(L_v\Omega_0)
=(v\nabla\omega+\omega  \div v)\Omega_0= \div (v\omega)\Omega_0\,.
$$

\bigskip
{\bf Appendix 3.}

\medskip
Here we derive (\der) and (\efi), for the convenience of the reader. Formulae (\der) are classical, see for example [A]. Inequality (\efi) is hardly new, but we were unable to find it in the literature. We recall the definitions

$$
F(\kappa)=\int_0^{\pi/2}{d\varphi\over{\sqrt{1+\kappa\sin^2\varphi}}},\quad E(\kappa)=\int_0^{\pi/2}\sqrt{1+\kappa\sin^2\varphi}\,d\varphi\,.
\eqno\hbox{\rm (\EF)}
$$
The calculation of $E'$  is straightforward:
$$
E'(\kappa)={d\over d\kappa}\intzpp\sqrt{1+\kappa\sin^2\vf}\,d\vf=\intzpp {\sin^2\vf\, d\vf\over{2\sqrt{1+\kappa\sin^2\vf}}}={1\over{2\kappa}}\intzpp{1+\kappa\sin^2\vf-1\over\sqrt{1+\kappa\sin^2\vf}}\,d\vf={1\over2\kappa}(E-F)\,.
$$

For $F'$ we obtain
$$
F'(\kappa)=\intzpp-{\sin^2\vf\,d\vf\over 2(1+\kappa\sin^2\vf)^{3\over2}}=\intzpp{{-1-\kappa\sin^2\vf+1}\over{2\kappa(1+\kappa\sin^2\vf)^{3\over2}}}\,d\vf=
-{F\over 2\kappa}+{1\over 2\kappa}\intzpp{d\vf\over(1+\kappa\sin^2\vf)^{3\over 2}}\,d\vf\,.
$$
To evaluate the last integral in terms of $E, F$, we set
$$
{1+\kappa\over 1+\kappa \sin^2\vf}=1+\kappa\sin^2 t\,\,.
$$
We calculate
$$
{d\vf\over(1+\kappa\sin^2\vf)^{3\over 2}}=-{1\over 1+\kappa}\sqrt{1+\kappa\sin^2 t}\,dt\,\,,
$$
which gives
$$
\intzpp{d\vf\over(1+\kappa\sin^2\vf)^{3\over 2}}\,d\vf={E\over 1+\kappa}\,.
$$
Hence
$$
F'={1\over 2\kappa}({E\over{1+\kappa}}-F)\,\,,
$$
as claimed. (One can also prove the identity by comparing the power series in $\kappa$.)

\medskip
Let us now turn to the proof of inequality (\efi). We first note that
$$
{d\over d\kappa}\,\,{E\over\sqrt{2+k}}= {E'\over{\sqrt{2+\kappa}}}-{E\over{2(2+\kappa)^{3\over 2}}}={1\over\sqrt{2+\kappa}}\left({E-F\over 2\kappa}-{E\over{2{(2+\kappa)}}}\right)={1\over{\kappa(2+\kappa)^{3\over 2}}}\left(E-(1+{\kappa\over 2})F\right)\,.
$$
Hence (\efi) is equivalent to showing that
$$
{d\over d\kappa}\left({E\over\sqrt{2+\kappa}}\right)<0,\qquad \kappa >0\,\,.\eqno\hbox{\rm (\dec)}
$$

In the second integral (\EF) which defines $E$ we can write $\sin^2\vf={1-\cos2\vf\over 2}$ and set $2\vf=\theta$ to obtain
$$
{E\over\sqrt{2+\kappa}}={1\over2\sqrt{2}}\int_0^\pi\sqrt{2+\kappa-\kappa\cos\theta\over 2+\kappa}\,d\theta={1\over2\sqrt{2}}\intzpp\left(\sqrt{1-\tau\sin\theta}+\sqrt{1+\tau\sin\theta}\right)\,d\theta\,,
$$
where
$$
\tau={\kappa\over{2+\kappa}}\,.
$$
We have
$$
{d\over d\tau} \intzpp\left(\sqrt{1-\tau\sin\theta}+\sqrt{1+\tau\sin\theta}\right)\,d\theta=
\pul\intzpp\left(-{\sin\theta\over\sqrt{1-\tau\sin\theta}}+{\sin\theta\over\sqrt{1+\tau\sin\theta}}\right)\,d\theta\,\,.
$$
As the last integral is obviously strictly negative for $0<\tau<1$ and $\tau$ is strictly increasing in $\kappa$, we have established (\dec), and hence (\efi) is proved.

\bigskip
\bigskip
\centerline{\it Acknowledgement} \smallskip The research was
supported in part by a grant from the National Science Foundation.

\parindent=1.5truecm
\bigskip
\bigskip
\centerline{REFERENCES}
\medskip

\item{[A]} Akhiezer, N. I.,\mez Elements of the theory of elliptic functions, {\it Translation of Mathematical Monographs, volume 79}, American Mathematical Society, 1990.

\item{[Am]} Amick, C.,\mez On Leray's problem of steady Navier-Stokes flow past a body in the plane, {\it Acta Math.} 161 (1988), no. 1--2, 71--130.

\item{[B]} Batchelor, G. K. \mez {\it An Introduction to Fluid
Dynamics}, Cambridge University Press, 1974 paperback edition.

\item{[CK]} Cannone, M. and Karch, G. \mez
           Smooth or singular solutions to the Navier-Stokes system?
           {\it J. Differential Equations} {\bf 197} (2004), no. 2, 247--274.

\item{[CY]} Chang, S.Y.A. and Yang, P.C.\mez The inequality of
Moser and Trudinger and applications to conformal geometry, {\it
Communications on Pure and Applied Mathematics,} Vol. LVI,
1135--1150, 2003.

\item{[DFN]} Dubrovin, B.A.,  Fomenko, A.T. and Novikov, S.P. {\it
Modern Geometry -- methods and applications}, Springer 1984--1990.

\item{[F]}
Finn, R., On the exterior stationary problem for the Navier-Stokes equations,
and associated perturbation problems,
 Arch.\ Rational Mech.\ Anal.\  19  1965 363--406.

\item{[FR]} Frehse, J, and R\accent23 u\v zi\v cka, M. \mez
Existence of regular solutions to the steady Navier-Stokes
equations in bounded six-dimensional domains, {\it Ann. Scuola
Norm. Sup. Pisa Cl. Sci.} (4) 23 (1996), no. 4, 701--719 (1997).

\item{[G]} Galdi, G.P. \mez {\it Introduction to the Mathematical
Theory of the Navier-Stokes Equations}, Volume II, Springer 1994.

\item{[Ha]} Hamel, G., \mez Spiralf\"ormige Bewegungen z\"aher Fl\"ussigkeiten, {\it  Jahresbericht d.\ Deutschen  Mathem.\-Vereinigung}, XXV, 34-65 (1917).

\item{[KS]} Korolev, A., Sverak, V.,\mez
On the large-distance asymptotics of steady state solutions of the Navier-Stokes equations in 3D exterior domains,
{\it Ann. Inst. H. Poincaré Anal. Non Linéaire} 28 (2011), no. 2, 303--313.

\item{[L]} Landau, L.D. \mez A new exact solution of the
Navier-Stokes equations, {\it Dokl.\ Akad.\ Nauk SSSR}, {\bf 43},
299, 1944.

\item{[LL]} Landau, L.D. and Lifschitz, E.M. \mez {\it Fluid Mechanics,
second edition,} Butterworth-Heinemann, 2000 paperback reprinting.

\item{[Le]} Leray, J. \mez Etude de Diverses \'Equations
Int\'egrales non Lin\'eaires et de Quelques Probl\`emes que Pose
l' Hydrodynamique, {\it J.\ Math.\ Pures Appl.}, {\bf 12}, 1-82,
1933.

\item{[NP]}
Nazarov, S. A., Pileckas, K. \mez On steady Stokes and Navier-Stokes problems with zero velocity at infinity
                               in a three-dimensional exterior domain, {\it J. Math.\ Kyoto Univ.} 40-3 (2000), 475--492.

\item{[Se]} Serrin, J. \mez The swirling vortex, {\it
Philosophical Transactions of the Royal Society of London,} Volume
271, 325--360, 1972.

\item{[Sq]} Squire, H.B. \mez The round laminar jet, {\it Quart.
J. Mech. Appl. Math.} 4, (1951). 321--329.

\item{[St1]} Struwe, M. \mez\mez Regular solutions of the stationary Navier-Stokes equations on $R^5$. {\it Math.\ Ann.} 302 (1995), no. 4, 719–741.

\item{[St2]} Struwe, M. \mez\mez personal communication, 1997.

\item{[Sv1]} \v Sver\'ak, V.  \mez\mez On Landau's Solutions of the Navier-Stokes Equations, arXiv:math/0604550, (2006).

\item{[Sv2]} \v Sver\'ak, V. \mez \mez unpublished note

\item{[T]} Tsai, T.P. \mez \mez Thesis, University of Minnesota,
1998.

\item{[TX]} Tian, G, and Xin, Z. \mez\mez One-point singular solutions
to the Navier-Stokes equations, {\it Topol. Methods Nonlinear
Anal.} 11 (1998), no. 1, 135--145.

\bigskip
\bigskip

 \end